\documentclass{gtart}

\def\ifplaintex{\expandafter\ifx\csname documentclass\endcsname\relax}

\def\gtp{{\mathsurround=0pt\it $\cal G\mskip-2mu$eometry \&\ 
$\cal T\!\!$opology $\cal P\!$ublications}}  

\def\recd{{\small Received:\qua\receiveddate\ifx\reviseddate\relax
\else\qquad Revised:\qua\reviseddate\fi\par}} 

\def\lognumber#1{\def\thelognumber{#1}}
\def\volumenumber#1{\def\thevolumenumber{#1}}
\def\volumeyear#1{\def\thevolumeyear{#1}}
\def\papernumber#1{\def\thepapernumber{#1}}
\def\pagenumbers#1#2{\def\startpage{#1}\def\finishpage{#2}}
\def\published#1{\def\publishdate{#1}}

\def\received#1{\def\receiveddate{#1}}

\def\accepted#1{\def\accepteddate{#1}}

\long\def\asciiabstract#1{\long\def\theasciiabstract{#1}}

\let\\\par\let\thelognumber\relax\let\thevolumenumber\relax
\let\thepapernumber\relax\let\thevolumeyear\relax\let\startpage\relax
\let\finishpage\relax\let\publishdate\relax\let\receiveddate\relax
\let\reviseddate\relax\let\accepteddate\relax\let\theasciititle\relax
\let\theasciiauthors\relax
\let\theasciiabstract\relax

\let\theasciiemail\relax

\ifplaintex
\font\logobig=cmssbx10 scaled 3836
\font\logomed=cmssbx10 scaled 2557
\else
\font\logobig=cmssbx10 scaled 4200
\font\logomed=cmssbx10 scaled 2800
\fi

\long\def\makeagttitle{   
\count0=\startpage
\agt\hfill      
\hbox to 45truept{\vbox to 0pt{\vglue -13truept{\logomed A\kern -.37em{\logobig 
T}\kern -.38em G}\vss}\hss}
\break
{\small Volume \thevolumenumber\ (\thevolumeyear)
\startpage--\finishpage\nl
Published: \publishdate}

\vglue .25truein

{\parskip=0pt\leftskip 0pt plus
1fil\def\\{\par\smallskip}{\Large\bf\thetitle}\par\medskip} \vglue
0.05truein

{\parskip=0pt\leftskip 0pt plus 1fil\def\\{\par}{\sc\theauthors}
\par\medskip}%
 
\vglue 0.03truein 

{\small\leftskip 25truept\rightskip 25truept{\bf Abstract}\stdspace\theabstract

{\bf AMS Classification}\stdspace\theprimaryclass
\ifx\thesecondaryclass\relax\else; \thesecondaryclass\fi\par
{\bf Keywords}\stdspace \thekeywords\par}\vglue 7truept

}   

\ifplaintex
\font\phead=cmsl9 scaled 950
\font\pnum=cmbx10 scaled 913
\font\pfoot=cmsl9 scaled 950
\headline{\vbox to 0pt{\vskip -4.5mm\line{\small\phead\ifnum
\count0=\startpage ISSN 1472-2739 (on-line) 1472-2747 (printed)
\hfill {\pnum\folio}\else\ifodd\count0\def\\{ }%
\ifx\theshorttitle\relax\thetitle\else\theshorttitle\fi\hfill{\pnum\folio}
\else\def\\{ and }{\pnum\folio}\hfill\ifx\theshortauthors\relax\theauthors
\else\theshortauthors\fi\fi\fi}\vss}}
\footline{\vbox to 0pt{\vglue 0mm\line{\small\pfoot\ifnum\count0=\startpage
\copyright\ \gtp\hfill\else
\agt, Volume \thevolumenumber\ (\thevolumeyear)\hfill\fi}\vss}}
\else
\font=cmsl9 scaled 1050
\font\lnum=cmbx10 
\font=cmsl9 scaled 1050
\makeatletter
\def\@oddhead{{\small\lhead\ifnum\count0=\startpage ISSN 1472-2739 
(on-line) 1472-2747 (printed)\hfill {\lnum\number\count0}\else\ifodd\count0
\def\\{ }\ifx\theshorttitle\relax \thetitle \else\theshorttitle\fi\hfill
{\lnum\number\count0}\else\def\\{ and }{\lnum\number\count0}
\hfill\ifx\theshortauthors\relax 
\theauthors\else\theshortauthors\fi\fi\fi}}\def\@evenhead{\@oddhead}
\def\@oddfoot{\small\lfoot\ifnum\count0=\startpage\copyright\ \gtp\hfill\else
\agt, Volume \thevolumenumber\ (\thevolumeyear)\hfill\fi}
\def\@evenfoot{\@oddfoot}
\makeatother
\fi
\let\maketitlepage\makeagttitle
\let\makeshorttitle\maketitlepage
\let\maketitle\maketitlepage

\newwrite\gtoutfile
\long\gdef\makeheadfile{  
{\def\\{, }\def\s{ }
\immediate\openout\gtoutfile head.xxx
\immediate\write\gtoutfile{To: math@arxiv.org}
\immediate\write\gtoutfile{Subject: put OR rep NNNNN:ppppp}
\immediate\write\gtoutfile{--text follows this line--}
\immediate\write\gtoutfile{Proxy-for: \ifx\theasciiauthors\relax
\theauthors\else\theasciiauthors\fi\s<\ifx\theasciiemail\relax\theemail\else\theasciiemail\fi>}
\immediate\write\gtoutfile{\noexpand\\}
\immediate\write\gtoutfile{Authors: \ifx\theasciiauthors\relax
\theauthors\else\theasciiauthors\fi}
{\def\\{ }\immediate\write\gtoutfile{Title: \ifx\theasciititle\relax
\thetitle\else\theasciititle\fi}}
\immediate\write\gtoutfile{Subj-class: GT or SG, GR etc}
\immediate\write\gtoutfile{MSC-class: \theprimaryclass\ifx\thesecondaryclass\relax\else, \thesecondaryclass\fi}
\immediate\write\gtoutfile{Journal-ref: Algebr. Geom. Topol. \thevolumenumber\s
(\thevolumeyear) \startpage-\finishpage}
\immediate\write\gtoutfile{Comments: Published by Algebraic and
Geometric Topology at}
\immediate\write\gtoutfile{\s\s\s  http://www.maths.warwick.ac.uk/agt/AGTVol\thevolumenumber/agt-\thevolumenumber-\thepapernumber.abs.html}
\immediate\write\gtoutfile{\noexpand\\}
\immediate\write\gtoutfile{}
\ifx\theasciiabstract\relax
\immediate\write\gtoutfile{\theabstract}\else
\immediate\write\gtoutfile{\theasciiabstract}\fi
\immediate\write\gtoutfile{}
\immediate\write\gtoutfile{\noexpand\\}
\immediate\write\gtoutfile{}
\immediate\closeout\gtoutfile}}  

\def\maketitlepage{\makeagttitle\makeheadfile}
\let\makeshorttitle\maketitlepage
\let\maketitle\maketitlepage

\lognumber{33}
\volumenumber{3}
\volumeyear{2003}
\papernumber{33}
\published{5 October 2003}
\pagenumbers{993}{1004}
\received{21 Febuary 2003}
\accepted{23 September 2003}

\usepackage{amssymb}
\newtheorem{theorem}{Theorem}
\newtheorem{lemma}{Lemma}[section]
\theoremstyle{definition}
\newtheorem{definition}{Definition}
\newcommand{\eps}{{\varepsilon}}
\newcommand{\C}{{\mathbf C}}

\newcommand{\R}{{\mathbf R}}
\newcommand{\Z}{{\mathbf Z}}

\begin{document}

\title {On three-periodic trajectories of\\multi-dimensional dual billiards}
\author{Serge Tabachnikov}
\address{Department of Mathematics\\Pennsylvania State 
University\\University Park, PA 16802, USA}
\email {tabachni@math.psu.edu}
\url{www.math.psu.edu/tabachni}

\begin{abstract}
We consider the dual billiard map with respect to a smooth strictly 
convex closed
hypersurface in linear $2m$-dimensional symplectic space and prove 
that it has at least $2m$ distinct
3-periodic orbits.
\end{abstract}
\asciiabstract{We consider the dual billiard map with respect to a 
smooth strictly convex closed hypersurface in linear 2m-dimensional
symplectic space and prove that it has at least 2m distinct
3-periodic orbits.}



\maketitle

\section{Introduction and formulation of results}

The dual billiard map is an outer counterpart of the usual billiard 
ball map.  In the plane, the
dual (or outer) billiard map $T$ is defined as follows. Let $M 
\subset \R^2$ be a smooth closed
strictly convex curve,  and let $z$ be a point in its exterior. There are
two tangent lines to  $M$ through $z$; choose one of them (say, the 
right one from the view-point
of $z$) and define  $T(z)$ to be the  reflection of $z$ in the point 
of tangency. J. Moser  put
forward the study of  dual billiards in \cite{Mo1, Mo2}, and by now, 
there exists a substantial
literature on the subject: see 
\cite{Am,Bo,D-T,G-K,G-S,K,S-V,Ta1,Ta2,Ta3,Ta4,Ta5,Ta6,Ta7}. Exotic
as they may seem, dual billiards are close relatives of the usual, 
inner ones; on the sphere, the
two systems are isomorphic via the projective duality.

Most of the research done on dual billiards concerned the plane case. 
However dual billiards can be
defined in multi-dimensional setting as well, \cite{Ta1,Ta2,Ta3}. 
This definition is not widely
known and is not immediately obvious, so we reproduce it here. One 
wants to replace the curve $M$
by a hypersurface in Euclidean space, and then the problem is that 
one has too many tangent lines to
$M$ through a given point. This difficulty is resolved as follows. 
One takes   the linear symplectic
space $(\R^{2m},\omega)$ where $\omega = \sum dx_i \wedge d y_i,\ 
i=1,\dots, m$, as the ambient space.
Let $M^{2m-1} \subset \R^{2m}$ be a smooth hypersurface and $q \in 
M$. The restriction of
$\omega$ to the tangent hyperplane $T_qM$ has a 1-dimensional kernel 
$\xi(q)$, called the
characteristic line. An orientation of $M$ induces an orientation of 
$\xi(q)$. This direction is
used as the tangent line at $q$ to define the dual billiard map.

\begin{definition} \label{corr} Two points $z_1, z_2 \in \R^{2m}$ are 
in the dual billiard
relation with respect to a smooth hypersurface $M$ if the mid-point 
$q=(z_1 + z_2)/2$ lies on $M$
and the line $(z_1 z_2)$ coincides with the characteristic line $\xi(q)$.
\end{definition}

In the plane (and for a convex dual billiard curve),  this is the 
definition of the dual
billiard map, given above. One can prove that the dual billiard 
relation is a symplectic relation.
One can also prove that if $M$ is a smooth closed strictly convex 
hypersurface then, for every
point $z_1$ outside of $M$, there exist exactly two points $z_2$ such 
that $z_1$ and $z_2$ are in
the dual billiard relation with respect to $M$. Moreover, for one 
choice of the point $z_2$, the
orientation of the ray $z_1 z_2$ coincides with that of $\xi(q)$, and 
for another the orientation
is opposite (just as in  the plane). Thus one defines the dual 
billiard map $T$, a
symplectomorphism of the exterior $M$; see \cite{Ta1,Ta2,Ta3}.

The dual billiard map may be
considered as a discretization of the characteristic flow on $M$; the 
relation between the dual
billiard map and the characteristic flow is the same as between the 
usual billiard map and the
geodesic flow on the boundary of the billiard table, see discussions 
in \cite{Ta1, Ta2} and \cite
{Ta3}, section 4.5.

An important problem concerning a dynamical system is to estimate the 
number of its periodic
trajectories. There is a natural action of the dihedral group $D_n$ 
on the set of
$n$-periodic trajectories $(z_1,\dots, z_n),\ T z_i = z_{i+1}$ for 
$i=1,\dots,n$: one can
cyclically permute the points $z_i$ or change their cyclic order to 
the opposite. When counting
periodic trajectories, we always count the number of $D_n$-orbits.

For strictly convex smooth plane dual billiards the situation is the 
same as for
the usual billiards inside a strictly convex closed smooth plane 
curve: {\sl for every period
$n \geq 2$ and a every rotation number $1 \leq r \leq n/2$, coprime 
with $n$, there exist
at least two distinct $n$-periodic dual billiard trajectories with 
the rotation number $r$.}
For the usual billiards, this is Birkhoff's theorem \cite{Bi}, and 
the result extends to
area-preserving twist maps of an annulus, of which the plane dual 
billiard map is a particular
case.

Concerning  multi-dimensional dual billiards, the only known result 
so far was the following
theorem, see  \cite{Ta1,Ta2,Ta3}: {\sl for every smooth closed 
strictly convex hypersurface in
linear symplectic space and every odd prime $n$, the dual billiard 
map has an $n$-periodic
trajectory.} This is a rather weak estimate that can be considerably improved.

A comparison with the usual billiards is in order. Historically the 
first estimate below of the number
of periodic trajectories in multi-dimensional convex billiards is due 
to I. Babenko \cite{Ba} where
the 3-dimensional case is considered. This case is technically 
hardest, and \cite{Ba} contained an
error. Satisfactory lower bounds were recently obtained  for a 
generic smooth closed strictly
convex billiard hypersurface $M^{m-1} \subset \R^{m}$ by M. Farber 
and the present
author in \cite{F-T1, F-T2}: {\sl for an odd $n \geq 3$, the number of distinct
$n$-periodic billiard trajectories inside $M^{m-1}$ is not less than 
$(n-1)(m-1)$}. Related results, in
particular, on the number of non-closed billiard trajectories going 
from one point to another are
contained in \cite{Fa,Ha}.

Our goal in this paper is to estimate below the number of 3-periodic 
trajectories of a strictly
convex smooth dual billiard. Period 3 is smallest possible for the 
dual billiard map, and  as
such, it is analogous to period 2 for the usual billiard. 
Two-periodic billiard trajectories are
also called diameters or double normals: these are chords of a 
hypersurface $M^{m-1} \subset
\R^{m}$ that are perpendicular to $M$ at both end points. If 
$M^{m-1}$ is smooth and strictly convex
then there are at least $m$ such chords; this is one of the earliest 
applications of the
Morse-Lusternik-Schnirelman theory. An extension to immersed 
submanifolds, not necessarily of
codimension 1, is contained in \cite{T-W}, and further 
generalizations and improvements in
\cite{Pu1, Pu2} and \cite{Du}, where the case of generalized 
3-periodic trajectories is
considered.

Our main result is as follows.

\begin{theorem} \label{main} For every smooth closed strictly convex 
hypersurface $M^{2m-1}$ in
$\R^{2m}$ the number of distinct 3-periodic trajectories of the dual 
billiard map is not less than
$2m$. This lower bound is sharp: for every $m$ there exists a smooth 
closed strictly convex
hypersurface $M^{2m-1} \subset \R^{2m}$ for which this number is equal to $2m$.
\end{theorem}

It is an interesting problem to extend this result to greater 
periods, even, as well as odd, that is,
to find a dual billiard analog of \cite{F-T1, F-T2}, and to include 
not necessarily convex immersed
hypersurfaces, in the spirit of \cite{Du, Ha, Pu1, Pu2}.

\section{Proofs}

Our proof makes use of the critical point theory.  The first step  is 
to find a function whose
critical points correspond to periodic trajectories. For the usual 
billiards, this function is the
perimeter length on inscribed polygons; for dual billiards it is the 
symplectic area counted with
certain multiplicities.

Consider a smooth hypersurface $M^{2m-1} \subset \R^{2m}$. Let $n$ be 
an odd number and
$(z_1,\dots,z_n)$ an $n$-tuple of points in $\R^{2m}$. Set: $q_i = 
(z_i + z_{i+1})/2$
where the index $i = 1,\dots,n$ is understood cyclically. The 
following result is contained in
\cite{Ta1,Ta2,Ta3}.

\begin{lemma} \label{gen} The points $(z_1,\dots,z_n)$ form an 
$n$-periodic orbit of the dual
billiard relation with respect to $M$ if and only if the points 
$q_1,\dots,q_n$ lie on $M$ and the
$n$-tuple $(q_1,\dots,q_n)$ is a critical point of the function on $M 
\times \dots \times
M$ given by the formula:
\begin{eqnarray} \label{genfun}
F(q_1,\dots,q_n) = \sum_{1 \leq i < j \leq n} (-1)^{i+j} \omega (q_i, q_j).
\end{eqnarray}
\end{lemma}

\begin{proof}
  By Definition \ref{corr}, points $z_k$ and $z_{k+1}$ are in the dual 
billiard relation
if and only if $q_k \in M$ and the vector $z_{k+1} - z_k$ is parallel 
to $\xi(q_k)$.

It is convenient to introduce the complex structure $J$:
$$\omega(u,v) = Ju \cdot v = -u \cdot Jv.$$
Let $\nu(q)$ be the unit normal vector to $M$ at point $q \in M$. 
Then the vector $J \nu(q)$ is
parallel to the characteristic line $\xi(q)$.

Consider $\bar q = (q_1,\dots,q_n) \in M \times \dots \times M$. Then
\begin{eqnarray} \label{dif}
dF(\bar q) = \sum_{1 \leq i < j \leq n} (-1)^{i+j} (J q_i\ dq_j - J q_j\ dq_i).
\end{eqnarray}
By the Lagrange multipliers method, the point $\bar q$ is critical 
for $F$ if and only if
\begin{eqnarray} \label{lag}
dF(\bar q) = \sum_k \lambda_k \nu(q_k)\ dq_k
\end{eqnarray}
where $\lambda_k$ are Lagrange multipliers. Fix $k$, collect the 
terms in front of $dq_k$ in
(\ref{dif}) and (\ref{lag}) and equate:
\begin{eqnarray} \label{eq}
(-1)^k \sum_{j=0}^{n-2} (-1)^j q_{k+1+j} = -\lambda_k J \nu(q_k).
\end{eqnarray}
Since $n$ is odd one can express $z$ in terms of $q$:
$$
z_i = \sum_{j=0}^{n-1} (-1)^j q_{i+j},
$$
and hence
\begin{eqnarray} \label{ex}
z_{k+1} - z_k = 2 \sum_{j=0}^{n-2} (-1)^j q_{k+1+j}.
\end{eqnarray}
Compare (\ref{ex}) with (\ref{eq}) to conclude that $\bar q$ is a 
critical point if and only if
$z_{k+1} - z_k$ has the characteristic direction at $q_k$ for all 
$k=1,\dots,n$. This means that
the  points $z_k$ and $z_{k+1}$ are in the dual billiard relation.
\end{proof}

Note that the function $F$ does not depend on the choice of the 
origin. Note also that the dihedral
group $D_n$ acts on the arguments of $F$ by permutations, and $F$ is 
an odd function with respect
to this action: $F$ is invariant under cyclic permutations and 
changes sign if the order is
reversed. This property distinguishes $F$ from the perimeter length 
function in the usual
billiard set-up which is invariant under $D_n$. For $n=3$, the case 
under consideration
here, $F$ is (negative) the symplectic area of a triangle inscribed 
into the dual billiard
hypersurface $M$:
\begin{eqnarray} \label{fn}
F(q_1,q_2,q_3) = \omega(q_2,q_1) + \omega(q_3,q_2) + \omega(q_1,q_2).
\end{eqnarray}

Assume now that $M$ is strictly convex so that the dual billiard map 
$T$ is defined. Most of the
critical points of $F(q_1,\dots,q_n)$ do not correspond to 
$n$-periodic orbits of $T$ but rather to
backtracking ``fake" orbits $(\dots ,z_{i-1}, z_i, z_{i+1}, \dots)$ 
with $z_{i-1} = z_{i+1}$. This
backtracking occurs when $q_{i-1} = q_i$. Therefore, to estimate the 
number of $n$-periodic
orbits of $T$ one needs to consider the critical points of $F$ on the 
cyclic configuration
space
$$G(M,n) = \{(q_1,\dots,q_n) \in M \times \dots \times M|\ q_{i-1} 
\neq q_i\ \ {\rm for}\  \
i=1,\dots,n
\}$$
whose topology was studied in \cite{F-T1}. This critical set is acted 
upon by $D_n$, and each
$D_n$ orbit corresponds to a single $n$-periodic orbit of $T$.

 From now on, $n=3$, the dihedral group $D_3$ coincides with the 
symmetric group $S_3$, and the
function $F$ is skew-symmetric, see (\ref{fn}). In particular, if 
$q_i = q_j$ for some $i \neq j$
then $F(q_1,q_2,q_3) =0$. To prove Theorem \ref{main} we use 
Morse-Lusternik-Schnirelman theory. One
cannot immediately apply this theory to $G(M,3)$ because this space 
is open. A way around this
difficulty is as follows.

Choose a sufficiently small generic $\eps >0$ and consider the 
manifold with boundary
$$U = \{(q_1,q_2,q_3) \in M \times M \times M|\ F(q_1,q_2,q_3)  \leq 
-\eps \}.$$
Then $U \subset G(M,3)$. The group of cyclic permutations $\Z_3$ acts on $U$
freely; let $V=U/\Z_3$. The function $F$ descends to a function $\bar F$ on
$V$, therefore the number of 3-periodic trajectories of the dual billiard map
is not less than the number of critical points of $\bar F$. The 
negative gradient of
$\bar F$ has the inward direction along the boundary $\partial V$, and
the usual estimates of the Morse-Lusternik-Schnirelman theory apply. 
In particular, the number of
3-periodic orbits of $T$ is not less than the sum of Betti numbers of 
$V$ and greater than the
cohomological length of this space.

The critical points of $F$ in the set $F(q_1,q_2,q_3) \geq \eps$ are 
in one-to-one correspondence
with those in $U$, due to skew-symmetry of the function $F$. However 
one may lose critical points on
the zero level set, that is, 3-periodic orbits that lie in an 
isotropic affine 2-plane. In fact, no
such orbits exist, as the next lemma shows.

\begin{lemma} \label{zer} A non-degenerate triangle with zero 
symplectic area, inscribed into $M$,
is not critical for the symplectic area functional.
\end{lemma}

\begin{proof}
  Assume that $q_1, q_2, q_3$ is an inscribed triangle with zero 
symplectic area,
corresponding to a 3-periodic orbit of $T$. Then the vector $q_3 - 
q_2$ is parallel to the
characteristic direction $\xi(q_1)$.
Since the symplectic area vanishes, $\omega(q_3 - q_2, q_2 - q_1) = 
0$, and hence
the vector $q_2 - q_1$ lies in the symplectic orthogonal complement 
to $\xi(q_1)$, that is, the
tangent hyperplane $T_{q_1} M$. This contradicts strict convexity of $M$.
\end{proof}

We identify the homotopy equivalence class of $V$ in two steps. Let $W$ be the
set of orthonormal 2-frames $(e_1,e_2)$ in $\R^{2m}$ with 
$\omega(e_1,e_2) > 0$. One has a free
$\Z_3$ action on $W$ by rotating a frame in the plane spanned by 
$e_1$ and $e_2$ through $2\pi/3$ in
the direction from $e_1$ to $e_2$.

\begin{lemma} \label{eq1}
$U$ is $\Z_3$-equivariant homotopy equivalent to $W$.
\end{lemma}

\begin{proof}
Given a triangle $q_1 q_2 q_3 \in U$, parallel translate it in
such a way that its center of mass is at the origin. Then dilate the position
vectors of the vertices of the triangle to unit vectors. We now have three
  pairwise distinct unit vectors $z_1,z_2,z_3$ in $\R^{2m}$ that form a triangle
of positive symplectic area. Let $P$ be the 2-plane of this triangle.

Consider the function
$$\phi(z_1,z_2,z_3) = |z_1 - z_2| + |z_2 - z_3| + |z_3 - z_1|.$$
This is a smooth function on the cyclic configuration space 
$G(S^1,3)$ of  the unit circle in the
plane $P$. The critical points of $\phi$ are 3-periodic billiard 
trajectories inside the circle,
that is, equilateral triangles. The open manifold $G(S^1,3)$ is 
homotopically equivalent to the
manifold with boundary given by the condition
$$|z_1 - z_2|  |z_2 - z_3|  |z_3 - z_1| \geq \delta$$
where $\delta>0$ is generic and sufficiently small. The gradient of 
$\phi$ has the
inward direction along the boundary and defines a retraction of
$G(S^1,3)$ to the set of inscribed equilateral triangles (see, e.g., 
\cite{F-T1} for
details).  An equilateral triangle in the plane $P$ is determined by 
its first vertex
$e_1$, and since the restriction of $\omega$ on $P$ does not vanish, 
the choice of $e_1$ uniquely
determines $e_2$ so that $\omega(e_1,e_2) > 0$.

The above construction provides a $Z_3$-equivariant map from $U$ to 
$W$; this map gives the desired
homotopy equivalence.
\end{proof}

Let $\lambda = (-1+\sqrt{-3})/2$ be the  primitive cube root of 1. 
Consider the standard action of
$\Z_3$ on the sphere $S^{2m-1} \subset \C^m$ by multiplying every 
coordinate by $\lambda$.

\begin{lemma} \label{eq2}
$W$ is $\Z_3$-equivariant homotopy equivalent to $S^{2m-1}$.
\end{lemma}

\begin{proof}
We already identified $\R^{2m}$ with $\C^m$ by introducing the operator $J$.
We want to show that $W$ is $\Z_3$-equivariant homotopy equivalent to 
the set of orthonormal
2-frames $(e_1,e_2)$ with $e_2 = J e_1$.

Consider the function
$\psi(e_1,e_2) = \omega(e_1,e_2)$
on the Stiefel manifold of orthonormal 2-frames. Note that $\psi$ is
$\Z_3$-invariant.
We claim that the critical points of $\psi$ consist of two critical manifolds,
the complex and anticomplex frames $e_2 = \pm J e_1$.  Indeed, let
$(e_1,e_2)$ be a frame and $(u_1,u_2)$ its variation. Then
\begin{eqnarray}
e_1 \cdot u_1 = e_2 \cdot u_2 = e_1 \cdot u_2 + e_2 \cdot u_1 = 0. \label{var}
\end{eqnarray}
One has:
$$d\psi (u_1,u_2) = J e_1 \cdot u_2 - J e_2 \cdot u_1.$$
Therefore, by (\ref{var}), $(e_1,e_2)$ is a critical point if and 
only if the vector
$(-J e_2,J e_1)$ is a linear combination of the three vectors
$(e_1,0), (0,e_2), (e_2,e_1)$.
This is possible only when $e_2 = \pm J e_1$.

To complete the proof of the lemma, one chooses a generic $\delta >0$ 
and replaces $W$ by a
manifold with boundary $\omega(e_1,e_2) \geq \delta$. Then the 
gradient of $\psi$ retracts $W$ to
the set of complex frames $(e_1,e_2)$ with $e_2 = J e_1$, and this is 
the desired $Z_3$-equivariant
homotopy equivalence.
\end{proof}

It follows from Lemmas \ref{eq1} and \ref{eq2} that $V$ is homotopy 
equivalent to the
lens space $S^{2m-1}/\Z_3$. The cohomological length of the lens space with
coefficients in $\Z_3$ is $2m-1$. It follows from the 
Lusternik-Schnirelman theory that the function
$\bar F$ has at least $2m$ distinct critical points on $V$. This 
proves the lower bound of
Theorem \ref{main}.

We construct a hypersurface in $\R^{2m}$ for which the number of 
3-periodic dual billiard
trajectories equals $2m$ by a small perturbation of the unit sphere 
$S^{2m-1} \subset \R^{2m}$.
We describe $M$ in terms of its support function.

The support function of a smooth strictly convex closed hypersurface 
$M^{2m-1}$ is the function on
the unit sphere $h:S^{2m-1} \to \R$ whose value on $z \in S^{2m-1}$ 
is the maximum of the linear
function $\phi_z (x) = x \cdot z$ on $M$. Let $q(z) \in M$ be the 
point at which the maximum is
attained. Then
\begin{eqnarray} \label{sup}
q(z) = h(z) z + \nabla h(z)
\end{eqnarray}
where the position vector $z$ and the gradient $\nabla h(z)$ are 
considered as vectors in $\R^{2m}$.
The vector $z$ is the outward unit normal vector to $M$ at $q$, and 
the characteristic direction at
$q$ is given by the vector $J z$. A triple of distinct points $z_1, 
z_2, z_3 \in S^{2m-1}$
corresponds to a 3-periodic orbit of the dual billiard map with 
respect to $M$ if and only if there
exist real numbers $a_1,a_2,a_3$ such that
\begin{eqnarray} \label{cl}
q_i  + a_i J z_i = q_{i+1} - a_{i+1} J z_{i+1},\ \ i=1,2,3,
\end{eqnarray}
where $q_i$ are related to $z_i$ as in (\ref{sup}).

\begin{lemma} \label{perp}
Let $f$ be a smooth $\Z_3$-invariant function on $S^{2m-1}$ that has 
$N$ critical
$\Z_3$-orbits. Consider the hypersurface $M$ with the support 
function $h = 1 + \eps f$. Then, for
sufficiently small $\eps$, the number of 3-periodic orbits of the 
respective dual billiard map does
not exceed $N$.
\end{lemma}

\begin{proof}
Consider first the case when $h=1$, that is, $M$ is the unit sphere. 
Recall that
$\lambda$ denotes the primitive cube root of 1. Then, for every
$z \in S^{2m-1}$, the triple $z_i = \lambda^i z$ is a solution of 
(\ref{cl}) with all $a_i = \sqrt{3}$,
and every solution of (\ref{cl}) is of this form.

Next let $h = 1 + \eps f$. We consider the linearizations of 
(\ref{sup}) and (\ref{cl}) in $\eps$. The
number of solutions of the linearized system provides an upper bound 
on the number of genuine
solutions.

We are looking for a perturbed solution in the form
\begin{eqnarray} \label{apr}
z_i = \lambda^i (z+ \eps v_i),\ a_i = \sqrt{3} + \eps \alpha_i,\ \ i=1,2,3,
\end{eqnarray}
where $z \in S^{2m-1}$ and $v_1, v_2, v_3 \in T_z S^{2m-1}$ (there is 
redundancy in the
choice of vectors $v_i$, for example, one may assume that $v_3=0$). All
computations are made modulo $\eps^2$.
Denote the common value of $f(z), f(\lambda z)$ and $f(\lambda^2 z)$ 
by $c$. Likewise, one has:
$$\nabla f(z) = \lambda^{-1} \nabla f(\lambda z) = \lambda^{-2} 
\nabla f(\lambda^2 z),$$
and we denote this tangent vector at $z$ by $w$.

Substitute (\ref{apr}) to (\ref{sup}) to obtain:
$$q (z_i) = \lambda^i [ z + \eps (v_i + c z + w ) ],\ \ i=1,2,3.$$
Next, substitute to (\ref{cl}) and collect terms linear in $\eps$:
\begin{equation} \label{fin}
\frac{(c + \sqrt{-1} \alpha_i) z + w + (1 + \sqrt{-3}) v_i}{\lambda} =
(c - \sqrt{-1} \alpha_{i+1}) z + w +(1 - \sqrt{-3}) v_{i+1}
\end{equation}
for $i=1,2,3$.

Note that $\lambda (1 - \sqrt{-3}) = 1 + \sqrt{-3}$, hence the 
coefficients of $v_i$ and $v_{i+1}$ in
(\ref{fin}) are equal. Add the three equations (\ref{fin}) to obtain:
$$\left(3(1-\lambda)c + (1+\lambda) \sqrt{-1} (\alpha_1 + \alpha_2 + 
\alpha_3)\right) z = 3 (\lambda
-1) w.$$ Since $w$ is perpendicular to $z$, we conclude that $w=0$. 
Thus $z$ is a critical point of the
function
$f$. Note that once such a critical point $z$ is chosen,  one finds 
from (\ref{fin})
that all $\alpha_i = \sqrt{3} c$ and then  solves the linear system 
(\ref{fin}) for vectors $v_i$.

To summarize, each critical $\Z_3$-orbit of the function $f$ 
contributes a solution to (\ref{fin}),
that is, corresponds, in the linear approximation, to a 3-periodic 
orbit of the dual billiard map.
\end{proof}

To prove the second statement of Theorem \ref{main} we construct a 
$\Z_3$-invariant function on
$S^{2m-1}$ with $2m$ critical $\Z_3$-orbits. Let 
$x_1,\dots,x_m,y_1,\dots,y_m$ be Euclidean coordinates
in $\R^{2m}$ so that $x_i + \sqrt{-1} y_i,\ i=1,\dots,m$ are complex 
coordinates in $\C^m$.

\begin{lemma} \label{func}
Let $a_1,\dots,a_m$ be distinct reals. Then, for sufficiently small 
$\eps$, the function
\begin{eqnarray} \label{functf}
f(x_1,\dots,x_m,y_1,\dots,y_m) = \sum_{i=1}^m a_i 
\frac{x_i^2+y_i^2}{2} + \eps \sum_{i=1}^m
\frac{x_i^3 - 3x_i y_i^2}{3}
\end{eqnarray}
is $\Z_3$-invariant on $S^{2m-1}$ and has $2m$ critical $\Z_3$-orbits.
\end{lemma}

\begin{proof}
The first term in (\ref{functf}) is invariant under multiplication by 
complex numbers with
modulus 1. The second term consists of ${\rm Re} \left((x_i + 
\sqrt{-1} y_i)^3\right)$, hence it is
invariant under multiplication by $\lambda$. To find critical points 
of $f$ on the unit sphere, we use
the Lagrange multipliers method:
\begin{eqnarray} \label{syst}
a_i x_i + \eps(x_i^2 - y_i^2) = \eta x_i,\ a_i y_i - 2 \eps x_i y_i = 
\eta y_i,\ \ i=1,\dots,m,
\end{eqnarray}
where $\eta$ is a Lagrange multiplier. For each $i$, the system 
(\ref{syst}) has four solutions:
\begin{eqnarray} \label{four}
x_i=\frac{\eta-a_i}{\eps}, y_i=0;\ x_i = -\frac{\eta-a_i}{2\eps},
y_i=\pm \frac{\sqrt{3}(\eta-a_i)}{2\eps}
\end{eqnarray}
and $x_i=y_i=0$. The solutions (\ref{four}) form a $\Z_3$-orbit for 
which $x_i^2 + y_i^2 = (\eta -
a_i)^2/\eps^2$.

For  $I \subset \{1,\dots,m\}$  consisting of at least 2 elements, 
consider the function
$$g_I (\eta) = \sum_{i \in I} (\eta - a_i)^2,$$
and let
$$\delta = \min_{\eta} \min_{|I|\geq 2} g_I (\eta).$$
Since all $a_i$ are distinct, $\delta >0$. Let $\eps$ be so small 
that $\delta > \eps^2$.

We claim
that for every solution of the system (\ref{syst}) there exists a 
unique index $i \in \{1,\dots,m\}$
for which (\ref{four}) holds, and $x_j=y_j=0$ for $j \neq i$. Indeed, 
assume that
(\ref{four}) occurs for indices in a set $I$ with $|I| \geq 2$. Since 
the domain of the function is the
unit sphere, one has:
$$1 = \sum_{i=1}^m (x_i^2 + y_i^2) = \sum_{i \in I} \frac{(\eta - 
a_i)^2}{\eps^2} \geq
\frac{\delta}{\eps^2},$$
contradicting the choice of $\eps$.

Finally, there are $m$ choices of an index $i \in \{1,\dots,m\}$, and 
for each choice one has $(\eta -
a_i)^2 = \eps^2$. This gives two values of $\eta$ and two respective 
solutions (\ref{four}). Thus the
function $f$ has $2m$ critical $\Z_3$-orbits.
\end{proof}

\medskip

{\bf Acknowledgments}\qua I am grateful to M. Farber for numerous 
stimulating discussions and to
D. Handron for explaining his work on billiard trajectories.
Part of this work was done at Max-Planck-Institut in Bonn; it is a 
pleasure to acknowledge the
Institute's invariable hospitality. This work was partially supported
by BSF and NSF grants.

\Addresses\recd

\medskip


\begin{thebibliography}

\vskip 3 mm

\bibitem {Am} {\bf E. Amiran,} {\it Lazutkin coordinates and 
invariant curves for outer billiards.} J.
Math. Phys., 36 (1995), 1232--1241.

\bibitem {Ba} {\bf I. Babenko,} {\it Periodic trajectories of 
three-dimensional Birkhoff billiards.}
Math. USSR, Sbornik, 71 (1992),  1--13.

\bibitem {Bi} {\bf G. Birkhoff,} {\it Dynamical systems.} Amer. Math. 
Soc. Coll. Publ., 9, 1927.

\bibitem {Bo} {\bf Ph. Boyland,} {\it Dual billiards, twist maps and 
impact oscillators.} Nonlinearity,
9 (1996), 1411--1438.

\bibitem {D-T} {\bf F. Dogru, S. Tabachnikov,} {\it On polygonal dual 
billiard in the hyperbolic
plane.}  Reg. Chaotic Dynamics, 8 (2003), 67--82.

\bibitem {Du} {\bf F. Duzhin,} {\it Lower bounds for the number of 
closed billiard trajectories of
period 2 and 3 in manifolds embedded in Euclidean space.}  Int. Math. 
Res. Notes, 2003, no. 8,
425--449.

\bibitem {Fa} {\bf M. Farber,} {\it Topology of billiard problems.} 
I, II. Duke Math. J., 115
(2002), 559--585, 587--621.

\bibitem {F-T1} {\bf M. Farber, S. Tabachnikov,} {\it Topology of 
cyclic configuration spaces and
periodic  trajectories of multi-dimensional billiards.} Topology, 41 
(2002),  553--589.

\bibitem {F-T2} {\bf M. Farber, S. Tabachnikov,}  {\it Periodic 
trajectories in 3-dimensional convex
billiards.} Manuscripta Math., 108 (2002), 431--437.

\bibitem {G-K} {\bf E. Gutkin, A. Katok,} {\it Caustics for inner and 
outer billiards.} Comm. Math.
Phys., 173 (1995), 101--134.

\bibitem {G-S} {\bf E. Gutkin, N. Simanyi,} {\it Dual polygonal 
billiards and necklace  dynamics.}
Comm. Math. Phys., 143 (1991),  431--450.

\bibitem {Ha} {\bf D. Handron,} {\it Generalized billiard paths and 
Morse theory for manifolds with
corners.}  Topology Appl., 126 (2002),  83--118.

\bibitem {K} {\bf R. Kolodziej,} {\it The antibilliard outside a 
polygon.} Bull. Pol.  Acad. Sci., 37
(1989), 163--168.

\bibitem {Mo1} {\bf J. Moser,} {\it Stable and random motions in 
dynamical systems.}  Ann. of Math.
Stud., 77, 1973.

\bibitem {Mo2} {\bf J. Moser,} {\it Is the solar system stable?} 
Math. Intell., 1 (1978),  65--71.

\bibitem {Pu1} {\bf P. Pushkar',}  {\it Diameters of immersed 
manifolds and of wave fronts.} C. R.
Acad.  Sci., 326 (1998), 201--205.

\bibitem {Pu2} {\bf P. Pushkar',} {\it A generalization of Chekanov's 
theorem. Diameters of immersed
  manifolds and wave fronts.}   Proc. Steklov Inst. Math., 1998, no. 2 
(221), 279--295.

\bibitem {S-V} {\bf A. Shaidenko, F. Vivaldi,} {\it Global stability 
of a class of  discontinuous dual
billiards.} Comm. Math. Phys., 110 (1987), 625--640.

\bibitem {Ta1}  {\bf S. Tabachnikov,} {\it Outer billiards.} Russian 
Math. Surv., 48 (1993), No. 6,
81--109.

\bibitem {Ta2} {\bf S. Tabachnikov,} {\it On the dual billiard 
problem.} Adv. Math., 115 (1995),
221--249.

\bibitem {Ta3}  {\bf S. Tabachnikov,} {\it Billiards.} Soc. Math. France, 1995.

\bibitem {Ta4} {\bf S. Tabachnikov,} {\it Poncelet's theorem and dual 
billiards.} L'Enseign. Math.,
  39  (1993), 189--194.

\bibitem {Ta5} {\bf S. Tabachnikov,} {\it Commuting dual billiards.} 
Geom. Dedicata, 53 (1994), 57--68.

\bibitem {Ta6} {\bf S. Tabachnikov,} {\it Asymptotic dynamics of the 
dual billiard transformation.}
  J. Stat. Phys., 83 (1996), 27--38

\bibitem {Ta7} {\bf S. Tabachnikov,} {\it Dual billiards in the 
hyperbolic plane.}  Nonlinearity, 15
(2002), 1051--1072

\bibitem {T-W} {\bf F. Takens, J.  White,}  {\it Morse theory of 
double normals of immersions.}
Indiana Univ. Math. J., 21 (1971/1972), 11--17.

\end{thebibliography}
\end{document}